 \documentclass[a4paper,10pt]{article}                                    
\usepackage{algorithm, algorithmicx}
\usepackage{graphicx} 
\usepackage{amsmath} 
\usepackage{amssymb}  
\usepackage{amsfonts}
\usepackage{color}

\usepackage{tabularx}
\usepackage{makecell}
\usepackage{rotating}
\usepackage{diagbox}

\newcommand{\orth}{\mathsf{orth}}

\newcommand{\argmin}[1]{\underset{#1}{\mathrm{argmin\,}}}

\newcommand{\xmin}{x^{\star}}
\newcommand{\xtrue}{\widetilde{x}}
\newcommand{\xinter}{x^{\sharp}}
\newcommand{\dy}{\delta_y}
\newcommand{\dA}{\delta_A}
\newcommand{\DA}{\Delta_A}
\newcommand{\dAj}{\delta_{A,j}}
\newcommand{\dAi}{\delta_{A,i}}

\newcommand{\Dy}{\Delta_y}
\newcommand{\qA}{\overline{A}}
\newcommand{\qy}{\overline{y}}

\newcommand{\supp}{S}

\newcommand{\can}{\beta} 

\newcommand{\sign}{\text{sign}}

\newcommand{\R}{\mathbb{R}}
\newcommand{\N}{\mathbb{N}}

\newcommand{\one}{\mathbf{1}}

\renewcommand{\P}{\mathrm{n_a}}
\newcommand{\Q}{\mathrm{n_b}}

\newtheorem{lemma}{Lemma}

\newtheorem{theorem}{Theorem}

\newtheorem{remark}{\noindent \textbf{Remark}}{\normalfont}{\normalfont}
{\normalfont}{\normalfont}
 {\normalfont }{\normalfont}
{\normalfont}{\normalfont}

\newenvironment{proof}{{\bf Proof }}{\hfill $\Box$ \\}

\begin{document}

\title{Sparse linear regression from perturbed data} 
\author{S. M. Fosson, V. Cerone,  D. Regruto}
\maketitle

\begin{abstract}                         
The problem of sparse linear regression is relevant in the context of linear system identification from large datasets. When data are collected from real-world experiments, measurements are always affected by perturbations or low-precision representations. However, the problem of sparse linear regression from fully-perturbed data is scarcely studied in the literature, due to its mathematical complexity. In this paper, we show that, by assuming bounded perturbations, this problem can be tackled by solving low-complex $\ell_2$ and $\ell_1$ minimization problems. Both theoretical guarantees and numerical results are illustrated in the paper.
\end{abstract}
\section{Introduction}\label{sec:intro}
The problem of fully-perturbed linear regression, i.e., the solution of systems $Ax=y$ where $y\in\R^m$ and $A\in\R^{m,n}$ are both subject to errors, has long history and ubiquitous applications. Although less popular than the model where errors are confined to $y$, the fully-perturbed case is more realistic, since,  usually, also $A$ is either experimentally observed or transmitted from remote devices, thus prone to numerical and measurement inaccuracies. Nowadays, this fact is particularly noteworthy in data-driven modeling, where complex systems are described uniquely from collections of observed data. see \cite{brubook19} for an overview.

Fully-perturbed problems, also known as  errors-in-variables (EIV) problems, are encountered and studied in different areas, e.g., numerical analysis, signal processing, and system identification, see \cite[Section 2.5]{mar07} for a  list of examples. 
In system identification, which is the main motivating example for this work, the fully-perturbed problem arises in those linear static/dynamic systems with both input and output corrupted by noise, as illustrated in \cite{sod07}. The most popular approach is total least squares (TLS), proposed in \cite{gol73,gol80}, which  performs a maximum likelihood estimation on EIV regression. TLS is optimal under specific probabilistic assumptions on the perturbations, while, in absence of a probabilistic setting, the unknown but bounded (UBB) noise paradigm is usually considered. In particular, in \cite{pig12} a set membership identification scheme in presence UBB noise is developed via polynomial optimization \cite{las01,lasbook}. 

A drawback of the above mentioned methods is the need of a large number of samples $m>n$. Even though the storage of large amounts of data is not an issue, the processing has two main flaws: the complexity might prevent a real-time identification, e.g., in tracking time-varying systems, and over-parameterization may occur. Over-parameterization basically consists in the scarce capability of extrapolating the essential model from data and, possibly leading to a wrong physical description; this problem largely studied in machine learning and data-driven science,  see \cite{bru16,brubook19}. In linear system identification, this turns into the presence of redundant estimated parameters, which, although mathematically consistent, might be too complex and distant from the real physical behavior. In fact, as discussed in \cite{fat18}, in many cases, the essential behavior of large dimensional systems can be adequately described by models with a number of parameters significantly lower than the number of state variables.

As to linear systems, sparse models, i.e. models with few active parameters, can be obtained from a small number of samples $m<n$, as proven in compressed sensing (CS,  \cite{can06}). In a nutshell, CS states the conditions to recover a sparse vector $x\in\R^n$, i.e., a vector with $k\ll n$ non-zero components, from  linear regression $y=Ax+\eta\in\R^m$, with $m<n$, $\eta$ being a measurement noise; see \cite{fou13} for a complete overview. CS is also applied to system identification: in \cite{gu09,tot11,roj14,fat18}, different approaches, based on either $\ell_0$  or $\ell_1$ minimization, are proposed to identify sparse linear systems from few measurements.

In classical CS and sparse system identification, $A$ is assumed to be be exactly known. 
More recently, the EIV sparse linear regression problem is considered  \cite{fox19acc}. The EIV case is difficult as intrinsically non-convex. However, in \cite{fox19acc}, it is shown that the problem can be formulated as a linear program (LP) if the parameters are a priori known to be either non-positive or non-negative. Although  realistic in several applications, ranging from sensor selection to localization, see e.g., \cite{cal16,bay15}, this knowledge is not usually plausible in system identification.
\subsection{Contribution}
The goal of this work is to tackle the problem of sparse linear regression from fully-perturbed data and few measurements, by extending the contribution of \cite{fox19acc} to the general case of parameters with unknown signs. This makes the approach valuable for system identification. Specifically, we develop a low-complex pre-processing algorithm that provides information on the signs; the overall solution is then obtained by solving two convex problems. Sufficient theoretical guarantees are provided for the success of the recovery, and the effectiveness is verified through numerical experiments as well. 

This work extends the results in \cite{fox19acc} in three  directions. First, we develop a method to priorly estimate the signs, and we provide theoretical conditions that guarantee its success, see Theorem \ref{theor1} in Section \ref{sec:ta}. Second, we extend the analysis of robustness sketched in \cite[Section IV]{fox19acc} by providing a more complete discussion on sufficient conditions for success, intended as corrected selection of the relevant parameters, see Theorem \ref{prop1} in Section \ref{sec:ta}. Third, extended numerical experiments are illustrated.

The paper is organized as follows. In Section \ref{sec:ps}, we formally illustrate the problem. In Section \ref{sec:lp}, we develop the proposed recovery procedure, and in Section \ref{sec:ta}, we prove theoretical conditions for its success. Section \ref{sec:nr} is devoted to numerical simulations, discussions, and comparison to the state-of-the-art methods. Finally, we draw some conclusions in Section \ref{sec:con}.
\section{Problem statement}\label{sec:ps}
In this section, we present the problem that we aim to solve. 
Let  $\xtrue\in\R^n$ be $k$-sparse, with $k\ll n$. The subset of non-zero components of $\xtrue$ is called support, and denoted as $\supp$.
Our goal is to recover $\xtrue$ from
\begin{equation}\label{y=Ax}
y=A\xtrue,~~y\in\R^m,~~A\in\R^{m,n}~~m<n.
\end{equation}
However, we assume that both $y$ and $A$ are experimentally observed, then affected by noise. More precisely, we assume that the available data are 
\begin{equation}
\qy=y+\dy,~~ \qA=A+\dA
\end{equation}
where $\dy\in\R^m$ and $\dA\in\R^{m,n}$  are UBB perturbations, with known bounds, that is,
\begin{equation}
\|\dy\|_{\infty}\leq \Dy,~~~  \|\dA\|_{\infty}\leq \DA. 
\end{equation}
The UBB paradigm is exploited in many engineering problems, as possible alternative to probabilistic models: whenever the perturbation does not follow a specific distribution, it makes sense to assume a prior information on its maximum magnitude.
For example, the UBB description is consistent with low-precision or quantized data, which are subject to an unknown adjustment in a known range, see \cite{fox19acc} for details. This is usual in transmissions systems. In this case, the $\DA$ and $\Dy$ are exactly known if the properties of the transmitter are known.
Within the automatic control framework, the UBB paradigm is exploited and studied in  EIV system identification \cite{sod18}, state estimation \cite{ale08},  model predictive control \cite{bem99book,ric05}, direct data-driven control \cite{cer17}, and factor analysis \cite{cic19}. In these applications, the error bound can be estimated either from previous information on the model and on the measurement devices or from available training datasets.

In principle, TLS can be applied to find a solution; however, TLS requires $m>n$ and does not encourage sparse solutions, as discussed, e.g., in \cite{fat18}. Sparsity can be obtained by adding an $\ell_1$-norm regularizer. The case of perturbations affecting only $y$ originates the Lasso problem \cite{tib96}, which dates back to the 1990s, and is widely studied in  CS, see \cite{fou13}, and in system identification, see \cite{roj14}.  Lasso has two formulations, which read as follows, respectively: $\min_x \|Ax-y\|_2^2+\lambda\|x\|_1$, where $\lambda>0$ is a design parameter, and $\min_x \|x\|_1$ s. t. $\|Ax-y\|_2^2\leq \eta$, where $\eta$ is a known bound on the measurement noise. The second formulation is denoted also as Basis Pursuit Denoising (BPDN$_2$), and is equivalent to the first one for suitable values of $\lambda$ and $\eta$, see \cite[Theorem 3.2]{fou13}. It is worth mentioning that, in \cite{don06infty}, the slightly different formulation: $\min_x \|x\|_1$ s. t. $\|Ax-y\|_{\infty}\leq \eta$ is studied which is  denoted here as BPDN$_{\infty}$.

Lasso, BPDN$_2$, and BDPN$_{\infty}$ are convex only if $A$ is exactly known. In contrast, the approach developed in this work extends to the EIV framework. 
Given $\qy\in\R^m$, $\qA\in\R^{m,n}$, $m<n$, $\DA>0$, and $\Dy>0$, we consider the following problem, first proposed in \cite{fox19acc}:
\begin{equation}\label{theproblem1}
 \begin{split}
  \min_{x\in\R^n}\|x\|_1~\text{ s. t. } &y=A x,~A\in\R^{m,n},~m<n\\
  &\qy=y+\dy\\
  &\qA=A+\dA\\
  &\|\dy\|_{\infty}\leq \Dy\\
  &\|\dA\|_{\infty}\leq \DA.
 \end{split}
\end{equation}
This can be interpreted as a perturbed/quantized CS problem, as illustrated in \cite{fox19acc}. In general, in the CS setting, the EIV model is barely considered, since perturbations on $A$ yield non-convex problems and significantly limit the compression capability of the system. In particular, the TLS approach cannot be applied as the number of measurements is too small. The main works on EIV CS are \cite{her10,zhu11,yan12}. In \cite{her10}, the Lasso robustness to perturbations on $A$ is studied. In \cite{zhu11}, a sparse TLS approach is proposed, which adds a regularization term on $\dA$. In \cite{yan12}, a given structure on $\dA$ is established, i.e., $\dA=B\text{diag}(\beta_0)$, where  $B\in\R^{m,n}$ is known, and $\beta_0\in[-r,r]^n$, $r>0$, is unknown; thus, the direction of each column of $\dA$ is known, and the unknown variables are $n$ instead of $mn$. The related optimization problem \cite[Equation 11]{yan12} is based on Lasso, with an additive $\ell_{\infty}$-norm; it is biconvex and it is tackled through alternating minimization, which only achieves a local minimum. 

Beyond CS, Problem \eqref{theproblem1} envisages sparse system identification from compressed measurements of EIV static/dynamic linear systems.  The application to dynamic systems is now explained. Given the following SISO ARX model
\begin{equation}\label{ARX}
y_t=\sum_{p=1}^{\P} a_{p}y_{t-p}+\sum_{q=1}^{\Q} b_{q}u_{t-q}+e_t
\end{equation}
where $u_t$, $y_t$, and $e_t$ respectively are the input, the output, and the system noise, the aim is the estimation of  the parameters $\theta=(a_1,\dots,a_{\P},b_1,\dots, b_{\Q})^T$. If $m$ output observations are collected, the problem can be expressed as the solution to linear system $(y_{t+1},\dots, y_{t+m})^T=A\theta+(e_{t+1},\dots e_{t+m})^T$ where 
\begin{equation}\label{Toeplitz}
A=\left(\begin{array}{cccccc}
		y_{t}&\cdots&y_{t-\P+1}&u_{t}&\cdots&u_{t-\Q+1}\\
		y_{t+1}&\cdots&y_{t-\P+2}&u_{t+1}&\cdots&u_{t-\Q+2}\\
		\vdots&&&&&\vdots\\
		y_{t+m-1}&\cdots&y_{t+m-\P}&u_{t+m-1}&\cdots&u_{t+m-\Q}\\
	\end{array}\right).
\end{equation}
When $m<\P+\Q$, the problem is underdetermined; this is known in the literature as compressive system identification, see \cite{san11,tot11,car13}. In those papers, the problem is tackled by using Lasso and orthogonal matching pursuit approaches, and $A$ is assumed to be exactly known.

As to the EIV problem, in \cite{pig12}, the authors prove that the feasible set of a dynamic EIV system with bounded noise  can be tightly approximated using polynomial optimization, by  leveraging the specific Toeplitz structure of $A$. The static EIV approach, which does not leverage this structure, can be applied as well, while it provides a less tight feasible set, as shown in \cite{pig12}. In this work, we consider the static approach, while a dynamic approach that exploits the structure is left for future work.
\section{$\ell_2$ + $\ell_1$ approach}\label{sec:lp}
In this section, we present the proposed approach to the solution of Problem \eqref{theproblem1}.  If the signs are known, Problem \eqref{theproblem1}  can be relaxed to an LP problem, as illustrated in \cite{fox19acc}. To estimate the signs, we propose an $\ell_2$ regularized strategy. We denote by $\ell_2 + \ell_1$ the overall algorithm, that will be obtained in the end of this section.
In the following, we describe in detail the $\ell_1$ and the $\ell_2$ stages of the algorithm.
\subsection{Sparse linear regression via $\ell_1$ minimization for given signs}\label{sec:l1}
We now illustrate how to formulate the sparse linear regression problem when information on the signs of the true $\xtrue$ is available.
In \cite{fox19acc}, Problem \eqref{theproblem1} is rewritten as follows
\begin{equation}\label{sys2}
 \begin{split}
  \min_{x\in\R^n}\|x\|_1~\text{ s. t. } &\|\qy-(\qA-\dA)x\|_{\infty}\leq \Dy\\
  &\|\dA\|_{\infty}\leq \DA.
 \end{split}
\end{equation}
This formulation is easily obtained by observing that $\dy=\qy-y=\qy-Ax=\qy-(\qA-\dA)x$.
Let us define a generalized vector of signs $s(v)\in \{-1,1\}^n$ of a vector $v\in\R^n$ as follows: $s_i(v)=1$ if $v_i\geq 0$, and $s_i(v)=-1$ if  $v_i\leq 0$.
We call it ``generalized'' because, differently from the usual sign operator, it does not distinguish the null values; as a consequence, $s(v)$ is not unique for $v$. We notice that $v=s(v)\bullet |v|$, where $\bullet$ denotes the component-wise product, and $|v|=(|v_1|,\dots,|v_n|)^T.$

If we assume to know $s(\xtrue)$, we can replace $x$ with $s(\xtrue)|x|$ in \eqref{sys2}. Let $\preceq$ denote the component-wise inequality between matrices and $\one_n:=(1,\dots,1)^T\in\R^n$.  We have
\begin{equation}\label{relax}
\begin{split}
\|\qy-(\qA-\dA)x\|_{\infty}&\leq \Dy\\
&\Updownarrow\\
-\Dy\one_m\preceq \qy-(\qA-\dA)x&\preceq \Dy\one_m\\
&\Updownarrow\\
\left(\begin{array}{c}
         (\qA-\dA)x\\
         -(\qA-\dA)x\\
        \end{array}\right)&\preceq\left(\begin{array}{c}
         \qy+\Dy\one_m\\
         -\qy+\Dy\one_m\\
        \end{array}\right)\\
&\Downarrow\\        
\left(\begin{array}{c}
          \qA s(\xtrue)\bullet |x| -\DA \one_m\one_n^T|x|\\
          -\qA s(\xtrue)\bullet |x| -\DA \one_m\one_n^T|x|\\
        \end{array}\right)&\preceq\left(\begin{array}{c}
         \qy+\Dy\one_m\\
         -\qy+\Dy\one_m\\
        \end{array}\right)\\
&\Updownarrow\\        
\left(\begin{array}{c}
          \one_m s(\xtrue)^T\bullet\qA  |x| - \one_m\one_n^T\DA|x|\\
-\one_m \ s(\xtrue)^T\bullet\qA |x| -\one_m\one_n^T\DA|x|\\
        \end{array}\right)&\preceq\left(\begin{array}{c}
         \qy+\Dy\one_m\\
         -\qy+\Dy\one_m\\
        \end{array}\right).\\
        \end{split} 
\end{equation}
In the last system of equations, the unknown is $z=|x|$. Therefore, we can relax  Problem \eqref{sys2} to the following LP problem
\begin{equation}\label{p:lp}
 \begin{split}
  \min_{z\in\R_+^n}\sum_{i=1}^n z_i~&\text{ s. t. } C z \preceq g\\
  &\text{ where }\\
  &C= \left(\begin{array}{c}
          \one_m s(\xtrue)^T\bullet\qA-\DA\one_m\one_n^T\\
         - \one_m s(\xtrue)^T\bullet\qA-\DA\one_m\one_n^T\\
        \end{array}\right)\in\R^{2m,n}\\
    &g= \left(\begin{array}{c}
         \qy+\Dy\one_m\\
         -\qy+\Dy\one_m\\
        \end{array}\right)\in\R^{2m}.\\
        \end{split}
\end{equation}
As noticed in \cite{fox19acc}, in several applications $s(\xtrue)$ is known. For example, the case of non-negative parameters is often encountered in sparse optimization in image processing problems \cite{fra17}, power spectrum estimation \cite{voy16}, sensor selection \cite{cal16}, and localization \cite{bay15}. However, in a number of applications, e.g., in system identification, prior information on the signs is not available. When $s(\xtrue)$ is not known, the problem is not convex. More precisely, according to the results proposed in \cite{cer93}, in a single orthant of the parameter space,  the problem is a convex LP. Therefore, one could compute the overall solution by solving $2^n$ LP problems, one for each orthant. Nevertheless,  this is computationally intense for large $n$. It is then fundamental to develop a strategy to priorly estimate $s(\xtrue)$.
In this paper, we propose a strategy to estimate $s(\xtrue)$ based on $\ell_2$ regularization.

\subsection{Estimation of the signs via $\ell_2$ regularization}
\begin{figure}
\begin{center}
\includegraphics[width=0.7\columnwidth]{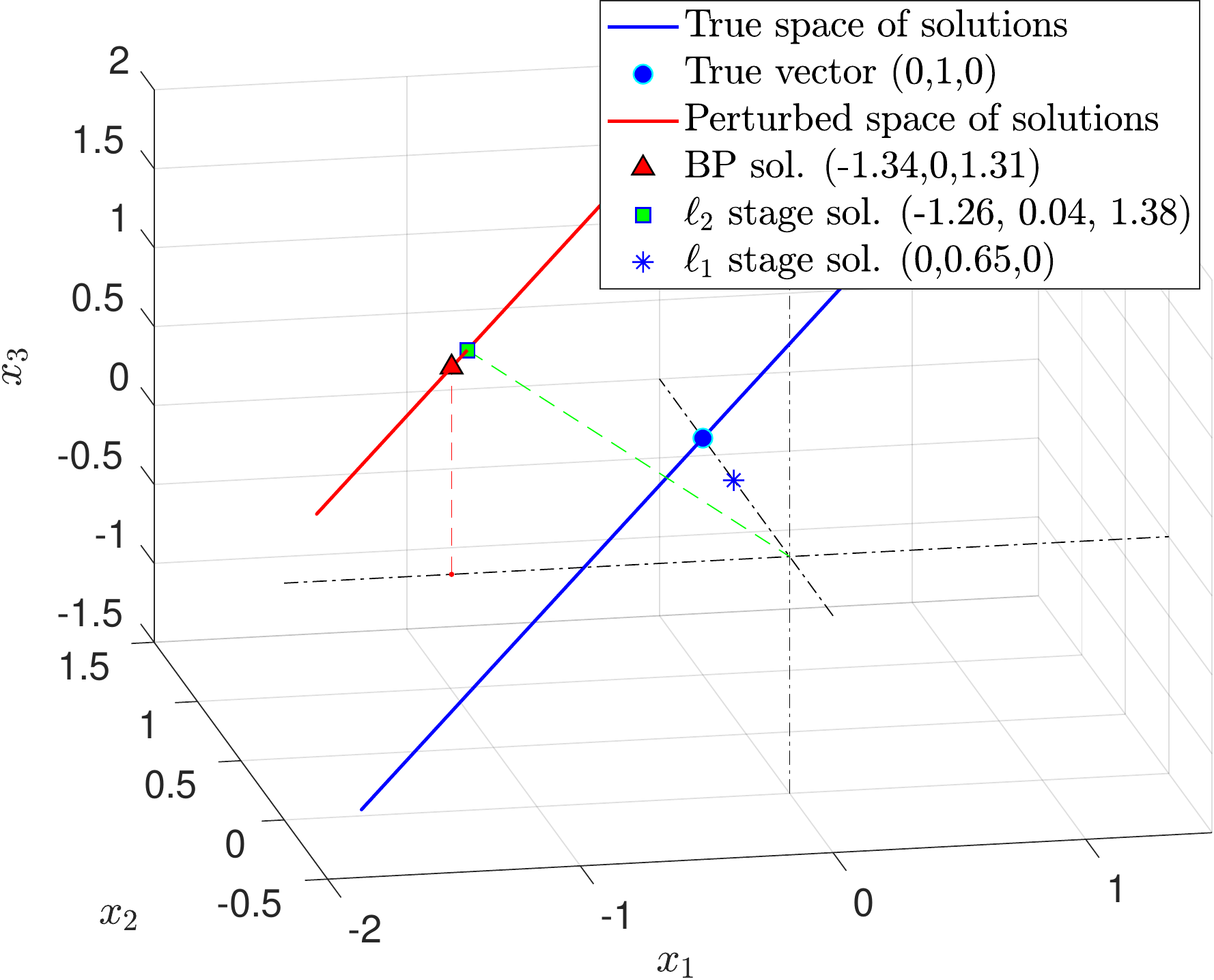}
\caption{Illustrating example: after perturbation, BP provides a wrong solution (red triangle). In contrast, the proposed $\ell_2$ stage correctly estimates the signs (green square). Given the signs, the proposed $\ell_1$ stage is run, which provides  the correct support (blue asterisk).}
\label{fig:motiv}
\end{center}
\end{figure}
In this section, we propose a Tikhonov ($\ell_2$-regularized) problem to estimate the signs. We introduce this approach with an example, while the rigorous analysis of its effectiveness is presented in Section \ref{sec:ta}. In particular, the example illustrates that perturbations modify the space of solutions in a way such that the sparsity pattern is compromised; this is a known issue in sparse optimization, CS, and factor analysis, see \cite{fou13,cic19}.

{\bf Example.} Let us consider $n=3$, $m=2$, $\xtrue=(0,1,0)^T$, and $$A=\left(\begin{array}{ccc}
                                                            1.5684&-2.5842&-0.1185\\-0.5477& 0.8012&0.1054
                                                           \end{array}\right).$$ Then, $y=( -2.5842,0.8012)^T$. The space of solutions of the equation $y=Ax$ is depicted in Fig. \ref{fig:motiv} as a blue line. Then, we consider perturbed data $\qy=(-2.4788,0.9580)^T$ and $$\qA=\left(\begin{array}{ccc}
                                                            1.7252&-2.8426&-0.1303\\-0.6025&0.8813&0.1159\\
                                                           \end{array}\right).$$  The space of solutions of $\qy=\qA x$ is the red line in Fig.  \ref{fig:motiv}. 
                                                        
Given $\qA$ and $\qy$, we can try to recover $\xtrue$ by using the popular Basis Pursuit (BP, \cite{fou13}), which consists in the convex problem: $\min_{x\in\R^n} \|x\|_1$ s. t. $\qA x=\qy$. 
Nevertheless, BP provides the solution $(-1.3378,0,1.3119)$, which misses the true support.
Then,  BP is not effective in this example. Therefore, we look for alternative strategies to   cope with perturbations. As already mentioned, the approach in \cite{fox19acc} is valuable for this purpose, but it requires information on the signs. Therefore, first of all, let us elaborate a strategy to get such information from the perturbed system.

In the example, we notice that $P =$ argmin $\|x\|_2$ s. t. $\qA x =\qy$ provides a right $s$ for $\xtrue$: specifically $P=(-1.2594,0.0444,1.3793)^T$, then   $s=(-1,1,1)^T$, which matches with $s(\xtrue)$. $P$, denoted as a green square in Fig. \ref{fig:motiv}, corresponds to the point of the red line closest to the origin. The intuition is the following: in presence of small perturbations, the space of solutions may substantially change its direction, but it still crosses the same orthants; thus, the $\ell_2$ minimization provides the right information on the signs.

Based on this intuition, we propose to estimate the signs by solving the $\ell_2$ problem: $$\min_{x\in\R^n}\|x\|_2~~~s.t.~~~ \qA x=\qy.$$
For practical purpose, we actually solve the relaxed problem
\begin{equation}\label{p:tik}
\min_{x\in\R^n} \|\qA x- \qy\|_2 + \lambda\|x\|_2^2
\end{equation}
with $\lambda>0$. If $\lambda$ is sufficiently small, the solution of Problem \eqref{p:tik} provides a good approximation of the constrained problem with the advantage that the solution can be computed in closed form.
This approach is rigorously analyzed in Section \ref{sec:ta}.

To conclude the example, given the correct signs, we apply the $\ell_1$ method as illustrated in Section \ref{sec:l1} (we assume $\DA=0.5168$ and $\Dy=0.3136$, which correspond to the double of the maximum perturbation. The obtained solution is (0,0.6445,0), which is a good approximation of $\xtrue$; in particular, it provides  the correct support.
\subsection{Summary of $\ell_2 +\ell_1$ approach}
After describing the two stages of the proposed approach, we summarize the $\ell_2 + \ell_1$ method in Algorithm \ref{alg}.

A rigorous theoretical analysis of Algorithm \ref{alg} is presented in Section \ref{sec:ta}, where sufficient conditions for success of the two stages are provided.
\begin{algorithm}
	\caption{\textbf{$\ell_2+\ell_1$} }\label{alg}
	\begin{algorithmic}[1] 
		\Statex Input: $\qy,\qA,\Dy,\DA$;
		\vspace{0.2cm}
		\Statex Output: $\xmin$ = estimate  of $\xtrue$
		\vspace{0.2cm}
		\State $\ell_2$ stage: solve Problem \eqref{p:tik} $~\Rightarrow~$ output: $\xinter$.
		\vskip0.15cm
		\State $\ell_1$ stage: solve Problem \eqref{p:lp} with $s(\xinter)$ instead of $s(\xtrue)$ $~\Rightarrow~$ output: $\xmin$
		%
	\end{algorithmic}
\end{algorithm}
\section{Analysis of $\ell_2+\ell_1$ approach}\label{sec:ta}
The aim of this section is to analyze the performance of the proposed algorithm $\ell_2+\ell_1$, summarized in Algorithm \ref{alg}. A key point of the two-stage procedure is the estimation of $s$ through Problem \eqref{p:tik}. Therefore, our first goal is to prove that Problem \eqref{p:tik} provides the correct signs, under suitable conditions. Afterwards, we elaborate on the performance of Problem \eqref{p:lp}.
%

In the following, we assume that the columns of $\qA$ are $\ell_2$-normalized. This assumption is not restrictive for sparse recovery: given any (non-normalized) $\qA$, let $D\in\R^{n,n}$ be the diagonal matrix with the $\ell_2$ norms of the columns of $\qA$ on the diagonal; then, $\qA D$ is normalized, and $\qA x = \qA D D^{-1} x$. Eventually, we  obtain a rescaled solution $D^{-1} x$. 

Moreover, we assume that the perturbed matrix $\qA$ has maximum rank $m$: this usually holds in the considered applications (in particular, in system identification), because $A$ is not expected to have a specific vector subspace structure, and by perturbing it is even less likely to achieve some algebraic structure.

Our analysis often refers to the concept of coherence of a matrix, which represents the normalized correlation among columns, see. e.g., \cite[Chapter 5]{fou13}. Given a matrix $M\in\R^{m,n}$, for any $m,n\in \N$, its coherence is defined as:
\begin{equation}\label{coherence}
\mu(M):= \max_{i,j\in\{1,\dots,n\}, i\neq j}\frac{|M_i^T M_j|}{\|M_i\|_2 \|M_j\|_2}
\end{equation}
where $M_i$ denotes the $i$-th column of $M$. Similarly, for any  subset $S\subset \{1,\dots,n \}$, $M_S$ is the submatrix obtained by selecting the columns in $S$.

%
\subsection{Analysis of $\ell_2$ stage}
Let us denote by $\orth(M)$ the operator that returns an orthogonal basis for the columns of $M$, see, e.g., \cite{fen09} for details.
The following theorem provides sufficient conditions for Problem \eqref{p:tik} to provide the right $s(\xtrue)$.  
%
\begin{theorem}\label{theor1}
Let $\xtrue \in \R^n$ be a $k$-sparse vector with support $S\subset \{1,\dots,n\}$. Let $c=\min_{i\in S}|\xtrue_i|$ and $d=\max_{i\in S}|\xtrue_i|$. Let $y=A\xtrue$, $y\in\R^m$, $A\in\R^{m,n}$, $m<n$, while the observed $\qy$ and $\qA$ are such that $\qy=y+\dy$, $\qA=A+\dA$, $\|\dy\|_{\infty}\leq \Dy$, $\|\dA\|_{\infty}\leq \DA$. $\qA$ is assumed to have maximum rank $m$ and $\ell_2$-normalized columns.

Let $Q:=\orth(\qA^T)^T$. If, for each $i\in S$, 
\begin{equation}\label{eq:theo1}
d\left|Q_i^T \sum_{j\in S\setminus \{i\}} Q_j\right|+ |f(\delta)|< Q_i^T Q_i c 
\end{equation}
where $f(\delta):=\qA^T(\qA \qA^T)^{-1} [\dy-\dA \xtrue]$, then the solution of Problem \eqref{p:tik} provides $s(\xtrue)$.
\end{theorem}
Before proving Theorem \ref{theor1}, we provide a technical lemma.
\begin{lemma}\label{lem1}
Given a matrix $A\in\R^{m,n}$, $m<n$, let us consider the singular value decomposition of $A^T=ULV^T$, $U\in\R^{n,n}$, $L\in\R^{n,m}$, $V\in\R^{m,m}$. Then, $L$ has $n-m$ null rows, and $UL=U_m L_m$ where $U_m\in\R^{n,m}$ is the submatrix of $U$ with columns corresponding the non-null rows of $L$, and $L_m\in\R^{m,m}$ is the diagonal matrix obtained by deleting the null rows. Let 
\begin{equation}
 Q:= U_m^T\in\R^{m,n}.
\end{equation}
Then,
\begin{equation}
A_S^T(AA^T)^{-1}A_S=Q_S^T Q_S. 
\end{equation}
\end{lemma}
\begin{proof}
It is known that $Q=\orth(A^T)^T$, see, e.g., \cite{fen09}. Then, in particular, $QQ^T=I$. 
As a consequence,
\begin{equation*}
\begin{split}
AA^T &= (Q^T L_m V^T)^T(Q^T L_m V^T)= VL_m^2 V^T  
\end{split}
\end{equation*}
and
\begin{equation*}
\begin{split}
A_S^T (A&A^T)^{-1} A_S = Q_S^T L_m V^T (V (L_m^2)^{-1} V^T) V L_m Q_S \\
&= Q_S^T L_m V^T (V (L_m^2)^{-1} V^T) V L_m Q_S = Q_S^T Q_S.
\end{split}
\end{equation*}
\end{proof}
Given Lemma \ref{lem1}, we can prove Theorem \ref{theor1}.

\begin{proof} (of Theorem \ref{theor1}). 
Let us consider the solution of the least squares problem with $\ell_2$ regularization:
\begin{equation}
z=\argmin{x\in\R^n} \|\qA x- \qy\|_2 + \lambda\|x\|_2^2
\end{equation}
where $\lambda>0$ is a sufficiently small design parameter. $z$ can be computed in closed form as:
\begin{equation}
z= (\qA^T\qA +\lambda I)^{-1}\qA^T \qy 
\end{equation}
where $I$ denotes the identity matrix of consistent dimensions. 
By exploiting the matrix inversion lemma, we can reduce the dimensionality of the involved matrix inversion:
\begin{equation}
(\qA^T\qA +\lambda I)^{-1}=\frac{1}{\lambda}\big[I- \qA^T (\qA \qA^T+\lambda I)^{-1}\qA\big].
\end{equation}
Then, we have:
\begin{equation}
\begin{split}
z&= \frac{1}{\lambda}\big[\qA^T \qy- \qA^T (\qA \qA^T+\lambda I)^{-1}\qA \qA^T \qy\big]\\
&= \frac{1}{\lambda}\big[\qA^T \qy- \qA^T (\qA \qA^T+\lambda I)^{-1}(\qA \qA^T\pm \lambda I) \qy\big]\\
&=\frac{1}{\lambda}\big[\qA^T \qy- \qA^T \qy -\qA^T(\qA \qA^T+\lambda I)^{-1}(-\lambda I)  \qy\big]\\
&=\qA^T(\qA \qA^T+\lambda I)^{-1} \qy.
\end{split}
\end{equation}

If $\qA$ has maximum row rank, $\qA \qA^T$ is non-singular, and all its eigenvalues are positive. Therefore, if $\lambda$ is sufficiently small, we can use the approximation: $(\qA \qA^T+\lambda I)^{-1} \approx (\qA \qA^T)^{-1}$. We remark that the same can not be said for $(\qA^T\qA+\lambda I)^{-1}$, since $\qA^T\qA$ is necessarily singular, $\qA^T\qA+\lambda I$ has eigenvalues equal to $\lambda$; therefore, a very small $\lambda$ affects the stability of the inversion.

Since $\qy=\qA \xtrue -\dA \xtrue +\dy$ and $\qA \xtrue = \qA_S \xtrue_S$, we then have
\begin{equation}
\begin{split}
 z&\approx \qA^T(\qA \qA^T)^{-1} \qy\\
 &= \qA^T(\qA \qA^T)^{-1} [\qA_S \xtrue_S -\dA \xtrue +\dy].
\end{split}
\end{equation}
By Lemma \ref{lem1}, $\qA_S^T(\qA \qA^T)^{-1} \qA_S=Q_S^T Q_S$, then
\begin{equation}
\begin{split}
z_S&\approx Q_S^T Q_S\xtrue_S+   \qA^T(\qA \qA^T)^{-1} [\dy-\dA \xtrue].
\end{split}
\end{equation}
Now, let us evaluate the distance between $z_S$ and $\xtrue_S$:
\begin{equation}\label{distanza}
\begin{split}
z_S-\xtrue_S& \approx (Q_S^T Q_S-I)\xtrue_S+   f(\delta).
\end{split}
\end{equation}
Our goal is to determine sufficient conditions so that, for each $i\in S$, $\sign(z_i)=\sign(\xtrue_i)$, which is equivalent to $|z_i-\xtrue_i|<|\xtrue_i|$.

By using \eqref{distanza}, the inequality that we have to solve, for $i\in S$, is:
\begin{equation}\label{distanza2}
\begin{split}
|Q_i^T Q_S\xtrue_S-\xtrue_i+ f(\delta)|< |\xtrue_i|
\end{split}
\end{equation}
where $Q_i^T$ is the $i$th column of $Q$.
By construction, the components of $Q_S^TQ_S$ are in $[-1,1]$, thus $1-Q_i^TQ_i\geq 0$. Now, by upper bounding the left-hand side of \eqref{distanza2}, we elaborate a sufficient condition to fulfill it:
\begin{equation}\label{distanza3}
\begin{split}
&\left|(Q_i^T Q_i -1)\xtrue_i+Q_i^T \sum_{j\in S\setminus \{i\}} Q_j \xtrue_j+ f(\delta)\right|< |\xtrue_i|\\
&(1-Q_i^T Q_i)|\xtrue_i|+\left|Q_i^T \sum_{j\in S\setminus \{i\}} Q_j \xtrue_j\right|+ |f(\delta)|< |\xtrue_i|\\
&\left|Q_i^T \sum_{j\in S\setminus \{i\}} Q_j \xtrue_j\right|+ |f(\delta)|< Q_i^T Q_i|\xtrue_i|\\
&d\left|Q_i^T \sum_{j\in S\setminus \{i\}} Q_j\right|+ |f(\delta)|< Q_i^T Q_i|\xtrue_i|.
\end{split}
\end{equation}
The thesis follows from the fact that for any $i\in S$, $|\xtrue_i|\geq c$.
\end{proof}
\begin{remark}\label{rem:0}
In \cite{fen09}, an orthogonalizing pre-processing is applied to work with $Q$ instead of $A$, based on the observation that the coherence of $Q$ is smaller than the coherence of $A$, as $Q$ is a tight frame. The research on tight frames is currently active, see, e.g., \cite{ela07,tsi14,dat19}, as they can  approximate equiangular frames, which are the ones that achieve the Welch bound, i.e. the minimal coherence, see \cite[Chapter 5]{fou13}.

From this observation, for $i,j\in S$,  we have $|Q_i^T Q_j|\leq \|Q_i\|_2 \|Q_j\|_2\mu(Q_S) \leq \|Q_i\|_2 \|Q_j\|_2 \mu(\qA_S)$, which can be used  to simplify \eqref{eq:theo1} as follows:
\begin{equation}\label{distanza4}
\begin{split}
&d(k-1)\|Q_j\|_2\mu(\qA_S)+ |f(\delta)|<  \|Q_i\|_2 c.
\end{split}
\end{equation}
%
%
\end{remark}
\begin{remark}\label{rem:1}
In practice, conditions \eqref{eq:theo1} or \eqref{distanza4} are not verifiable, as $\supp$ is not known. However, if some prior information on the structure of $A$ is available, one can experimentally estimate its coherence properties, then the sufficient number of measurements $m$ on average, see Remark \ref{rem:2} for a deeper discussion. On the other hand, also the range $[c,d]$ might be unknown a priori. However, it is reasonable to assume that one can estimate it from some physical knowledge of the parameters or from training data.
\end{remark}
\subsection{Analysis of $\ell_1$ stage}
Theorem \ref{theor1} provides conditions under which the first stage of $\ell_2+\ell_1$ is successful. Once the signs are correctly given, we can tackle the second stage, for which sufficient success conditions are provided in Theorem 1 in \cite{fox19acc}. Here, we propose a refined result, that states sufficient conditions to obtain the correct support. 
In the following, we assume that, in the final estimate, all the components with magnitude smaller than a given threshold $\tau$ are considered as null. This polishing is a common practice in sparse optimization, in particular when iterative solvers are used, and stopped before complete convergence. Such $\tau$ is usually assessed from a training dataset. In our case, if $c$ is known, it makes sense to set $\tau=\frac{c}{2}$ as threshold between zeros and non-zeros. Otherwise, we assume to be able to set a suitable lower bound $\tau<\frac{c}{2}$ from a training dataset.
\begin{theorem}\label{prop1}
Let us consider the setting of Theorem \ref{theor1}, and let us assume that Problem \eqref{sys2} provides the right $s(\xtrue)$. Let 
\begin{equation}
\phi:=\sqrt{m}(2\Dy+\DA (k+\widehat{k}) d) 
\end{equation}
where $\widehat{k}$ is the estimated sparsity level, and
\begin{equation}
\gamma(A):=\max_{i\in\supp}\sum_{l\in\supp,l\neq i} |A_i^TA_l|+\max_{j\in\supp^C}\sum_{l\in\supp}|A_j^TA_l|
\end{equation}
where $\supp^C$ is the complementary of $S$. Let $\tau$ be the above defined polishing threshold.
If there exists $\xi>0$ such that
\begin{equation}\label{GREATCONDITIONS}
\left\{\begin{aligned}
&\phi k \leq \tau\xi  \\
&\gamma(A)\leq 1-2\xi\\
\end{aligned}\right.
\end{equation}
then the solution of Problem \eqref{p:lp} exactly provides $S$, i.e., it selects the correct significant parameters.
\end{theorem}

\begin{figure*}
\begin{center}
\includegraphics[width=0.45\columnwidth]{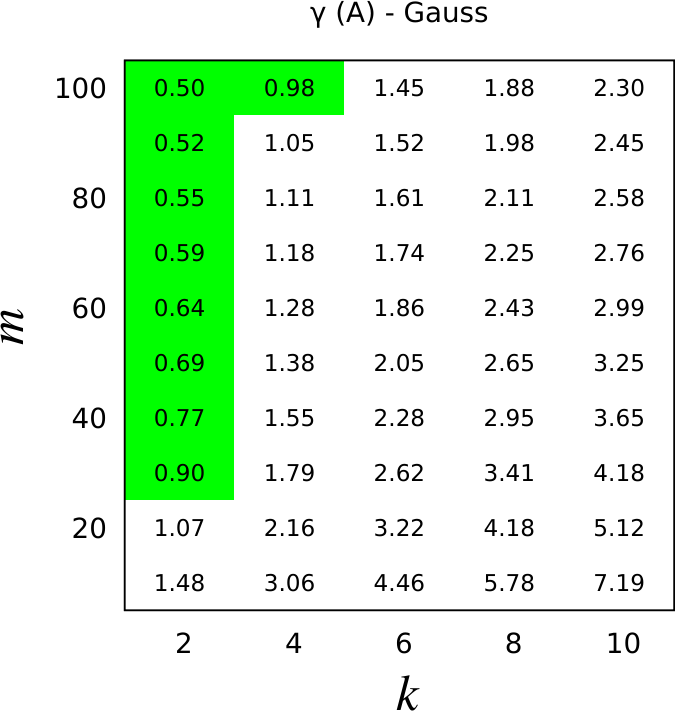}\hskip0.8cm
\includegraphics[width=0.45\columnwidth]{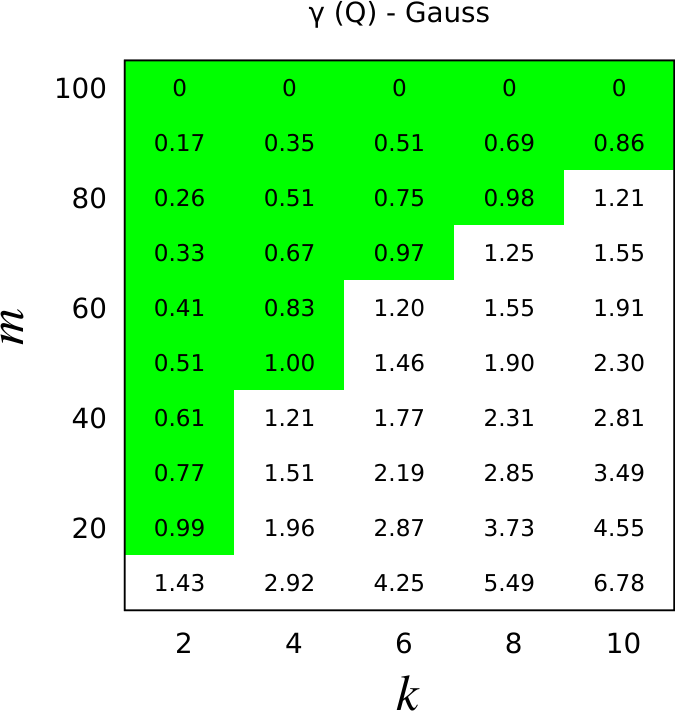}
\vskip1cm
\includegraphics[width=0.45\columnwidth]{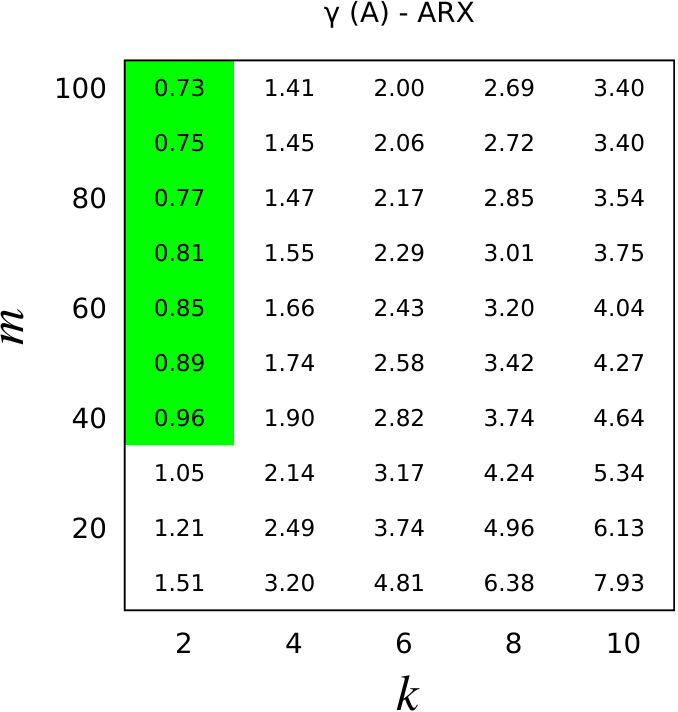}\hskip0.8cm
\includegraphics[width=0.45\columnwidth]{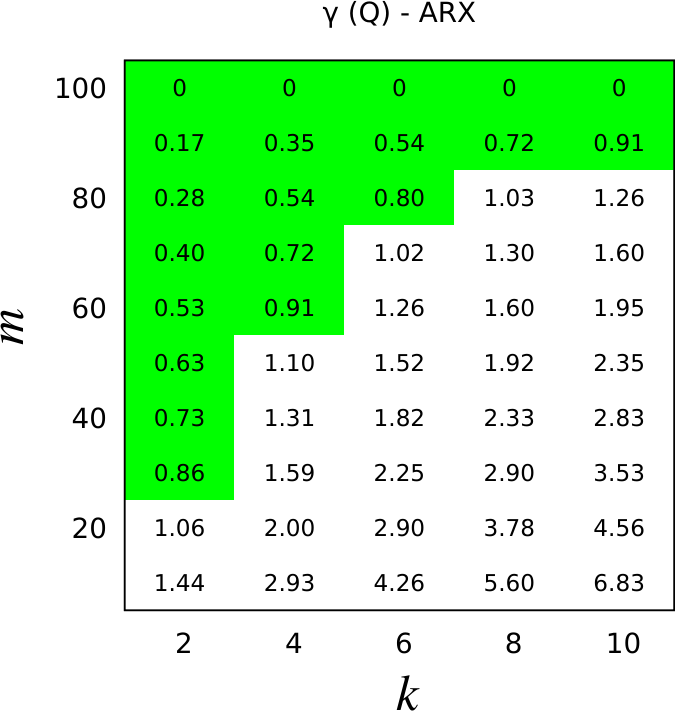}
\caption{$\gamma$ for $A\in\R^{m,n}$, $n=100$, and $Q\in\R^{m,n}$, obtained from $A$ through the orthogonalizing pre-processing by \cite{fen09}; $Q_S$. Standard Gaussian matrices (first row) and ARX matrices \eqref{Toeplitz} with Gaussian noise with variance $10^{-2}$ (second row) are considered. The results are averaged over 500 runs. The green color highlights the cases where $\gamma(A)<1$, which guarantees the exact support recovery for $2k\psi/c\leq \xi=(1-\gamma)/2$. As expected, 1) $Q$ is more successful than $A$, and 2) Gaussian, uncorrelated matrices are more successful than ARX matrices.}
\label{fig:coherences}
\end{center}
\end{figure*}

\begin{proof}
Let $\can\in\R^n$ be the solution of Problem \eqref{p:lp} and $w=\can-\xtrue$, where $\xtrue$ is the true vector to be recovered. 

If $\|w\|_{\infty}<\tau$ for a suitable threshold $\tau\leq \frac{c}{2}$, then the support of $\xtrue$ can be obtained from the support of $\can$. Therefore, our goal is to prove that $\|w\|_{\infty}<\tau$. We proceed by contradiction, by assuming that there is a component $j\in\{1,\dots,n\}$ such that
\begin{equation}\label{ass:absurd}
 |w_j|>\tau.
\end{equation}
The key idea of the proof is to formulate an LP problem in $|w|$, which shows that condition \eqref{ass:absurd} is inconsistent with the assumption that $\can$ solves Problem \eqref{p:lp}. In order to formulate this problem, we first notice that:
\begin{equation}\label{constraint2}
\begin{split}
\|\qA w\|_{\infty}&=\|\qA \can -\qA \xtrue\|_{\infty}\\
&\leq \|\qA\can-\qy\|_{\infty}+\|\qA\xtrue - \qy\|_{\infty}\\
&\leq\Dy+\DA\sum_{i=1}^n\can_i+\Dy+\DA\sum_{i=1}^n\xtrue_i\\
&\leq 2\Dy+\DA (k+\widehat{k}) d
\end{split}
\end{equation}
where $\widehat{k}$ is the sparsity level of $\can$, and can be upper bounded by $n$.

Then, given $|w|=(|w_1|\,\dots,|w_n|)^T$, we have
\begin{equation}\label{eq:splittino2}
\begin{split}
|w|&=|w + \qA^T \qA w- \qA^T \qA w| \preceq |\qA^T \qA w|+ |I-\qA^T \qA|~|w|.
\end{split}
\end{equation}
Now, let us define $$\Gamma:=I-|I-\qA^T \qA|\in\R^{n,n}.$$
It is worth noticing that, for $i\neq l$, $\Gamma_{i,l}=-|A_i^TA_l|<0$, while  $\Gamma_{i,i}=1-|1-\qA_i^T\qA_i|=1$.
Then, we can write
\begin{equation}\label{eq:pre_useful_inequality}
 \Gamma |w|\preceq  |\qA^T \qA w|.
\end{equation}
Furthermore, from \eqref{constraint2}, we obtain the following bound for $|\qA^T_i \qA w|$, for each $i\in\{1,\dots,n\}$
\begin{equation}\label{eq:useful_inequality}
\begin{split}
|\qA_i^T \qA w|&\leq \|\qA_i\|_2 \|\qA w\|_{2}\leq  \|\qA w\|_{\infty}\sqrt{m}\leq  \phi
\end{split}
\end{equation}
 From equations \eqref{eq:pre_useful_inequality} and \eqref{eq:useful_inequality}, we obtain
\begin{equation}\label{w_cond}
  \Gamma |w|\preceq \phi \one_n
\end{equation}
Moreover, as explained in \cite{don01,fox19acc}, the following holds:
\begin{equation}\label{use_w}
\|\can\|_1-\|\xtrue\|_1\geq\|w\|_1-2\|w_{\supp}\|_1  = \|w_{\supp^c}\|_1-\|w_{\supp}\|_1
\end{equation}
where $\supp$ is the support of $\xtrue$, and $\supp^c$ is its complementary.

As stated above, we assume that for some $j$, $|w_j|\geq \tau$. Then, by merging \eqref{w_cond} and \eqref{use_w}, we can formulate the following LP problem in $z:=|w|$:
\begin{equation}\label{z}
\begin{split}
\min_{z\in\R^n_{+}}& h^T z \\ 
&\text { s.t. } \left(\begin{array}{c}
                \Gamma\\
                -e_j^T\\
               \end{array}\right) z\preceq \left(\begin{array}{c}
                                              \phi \one_n\\
                                              -\tau
                                             \end{array}\right)
\end{split}
\end{equation}
where $h\in\{-1,1\}^n$ has components equal to $-1$ on $\supp$, and $+1$ elsewhere, so that $h^Tz=\|w_{\supp^c}\|_1-\|w_{\supp}\|_1$, and $e_j\in\R^n$ is $1$ in position $j$ and $0$ elsewhere.

If Problem \eqref{z} has a solution with positive penalty, then $\can$ is not a  solution of Problem \eqref{p:lp}, since from \eqref{use_w} its $\ell_1$ norm in not minimal in the feasible set. Then, to prove the contradiction, it is sufficient to show that a solution of  Problem \eqref{z} with positive penalty exists. To this purpose, we build the dual problem:
\begin{equation}\label{z2}
\begin{split}
\max_{\zeta\in\R^{n+1}_+}& (-\phi\one_n^T,~  \tau)\zeta\\
&\text { s.t. } \left(\begin{array}{cc}
                -\Gamma, &~e_j\\
               \end{array}\right)\zeta \preceq h\\
\end{split}
\end{equation}
For notational simplicity, let $\zeta=(u^T,\nu)^T$, $u\in\R_+^{n}$, $\nu\in\R_+$, so that:
\begin{equation}\label{uv}
\begin{split}
\max_{u\in\R^{n}_+,\nu\in\R_+}& -\phi \sum_{i=1}^n u_i+\tau\nu \\
&\text { s.t. } -\Gamma u + e_j\nu\preceq h\\
\end{split}
\end{equation}
Now, we show that there exists a choice of $u$ and $\nu$ such that, in Problem \eqref{uv},  the penalty is positive and  the constraints are satisfied. Based on  the zero duality gap between primal and dual in LP problems, see, e.g., \cite{lue16}, which implies that $\can$ is not a solution of Problem \eqref{p:lp}. 
Specifically, we propose a vector $u$ with support on $S$ and constant non-zero values: 
\begin{equation}\label{agoodsolution2}
u=\psi\one_S,~~\psi\in\R_+
\end{equation}
%
Then, Problem \eqref{uv} has positive penalty if
\begin{equation}\label{agoodsolution3}
-\phi k \psi+\tau\nu \geq 0
\end{equation}
along with the condition 
\begin{equation}\label{newlab}
-\Gamma u + e_j\nu=-\psi\Gamma_{\supp}\one_{\supp} + e_j\nu\preceq h. 
\end{equation}
We study this condition by distinguishing the cases $j \in \supp^C$ and $j\in \supp$.

{\textbf{Case 1:  $j \in \supp^C$}}

Condition \eqref{newlab} can be split as follows:
\begin{equation}\label{final_conditions1}
\begin{split}
&1)~\text{for each }i\in\supp:~~-\psi\sum_{l\in\supp}\Gamma_{i,l} \leq -1\\
&2)~\text{for each }i\in\supp^c,i\neq j:~~-\psi\sum_{l\in\supp}\Gamma_{i,l}\leq 1\\
&3)~j\in\supp^c:~~-\psi\sum_{l\in\supp}\Gamma_{j,l}+\nu \leq 1.
\end{split} 
\end{equation}

{\textbf{Case 2:  $j \in S$}}

The procedure is similar to that of Case 1, and yields
\begin{equation}\label{final_conditions2}
\begin{split}
&1)~\text{for each }i\in\supp, i\neq j:~~-\psi\sum_{l\in\supp}\Gamma_{i,l} \leq -1\\
&2)~j\in\supp:~~-\psi\sum_{l\in\supp}\Gamma_{j,l}+\nu \leq -1\\
&3)~\text{for each }i\in\supp^c:~~-\psi\sum_{l\in\supp}\Gamma_{j,l} \leq 1.
\end{split} 
\end{equation}
Finally, we merge Case 1 and Case 2. We notice that conditions \eqref{final_conditions1} and \eqref{final_conditions2} are both satisfied if
\begin{equation}\label{final_conditions3}
\begin{split}
&1)~\text{for each }i\in\supp^c:~~-\psi\sum_{l\in\supp}\Gamma_{i,l}+\nu \leq 1\\
&2)~\text{for each }i\in\supp:~~-\psi\sum_{l\in\supp}\Gamma_{i,l}+\nu \leq -1.
\end{split} 
\end{equation}
By recalling the definition of $\Gamma$ and the fact that its diagonal components are equal to 1, while the off-diagonal components are negative, the equations  \eqref{final_conditions3} are equivalent to
\begin{equation}\label{final_conditions4}
\begin{split}
&1)~\text{for each }i\in\supp^c:~~\sum_{l\in\supp}|A_i^TA_l| \leq \frac{1-\nu}{\psi}.\\
&2)~\text{for each }i\in\supp:~~\sum_{l\in\supp,l\neq i}|A_i^TA_l| \leq 1 -\frac{1+\nu}{\psi}
\end{split} 
\end{equation}
By considering the maximum of the left-hand sides of the two last inequalities, and by summing the two inequalities, we obtain \eqref{GREATCONDITIONS} with $\xi=\frac{\nu}{\psi}$.
Vice versa, it is easy to prove that \eqref{GREATCONDITIONS} imply \eqref{final_conditions4} if either $\psi=(\max_{i\in\supp^C}\sum_{l\in\supp}|A_i^TA_l|+\xi)^{-1}$ or $\psi=(\max_{i\in\supp}\sum_{l\in\supp,l\neq i}|A_i^TA_l|+\xi-1)^{-1}$, and $\nu=\xi\psi$.
\end{proof}
\begin{remark}\label{rem:2}
In practice, conditions \eqref{GREATCONDITIONS} are not verifiable, as $\supp$ is unknown. However, one can experimentally estimate their validity. As an example, in Fig. \ref{fig:coherences}, we evaluate $\gamma$ from data for Gaussian matrices and ARX systems \eqref{ARX}-\eqref{Toeplitz}.
\end{remark}
\begin{remark}\label{rem:3}
Similarly to Remark \ref{rem:0}, we can rewrite \eqref{GREATCONDITIONS} in terms of coherence of $\qA$. 
Concerning Case 1 in the proof of Theorem \eqref{prop1}, for each $i\in\supp$, $\sum_{l\in\supp}\Gamma_{i,l}\geq 1-(k-1)\mu(\qA_S)$. If  $1-(k-1)\mu(\qA_S)>0$, a sufficient condition to satisfy the first inequality in \eqref{final_conditions1} is 
\begin{equation}\label{cond_uno}
(k-1)\mu(\qA_S)\leq 1- \frac{1}{\psi}. 
\end{equation}
Similarly, a sufficient condition to satisfy the second inequality in \eqref{final_conditions1} is 
\begin{equation}\label{cond_due}
k\mu(\qA_{\supp^c})\leq\frac{1}{\psi}. 
\end{equation}

Finally, by setting $\nu<1$, a sufficient condition to satisfy the third inequality in \eqref{final_conditions1} is
\begin{equation}\label{cond_tre}
k\mu(\qA_{\supp^c})\leq \frac{1-\nu}{\psi}
\end{equation}
which is stronger than condition \eqref{cond_due}.
%
%
By the same procedure, for Case 2 we obtain the additional condition:
\begin{equation}\label{cond_quattro}
(k-1)\mu(\qA_S)\leq 1- \frac{1+\nu}{\psi}
\end{equation}
which is stronger than condition \eqref{cond_uno}.
By merging these observations, we conclude that conditions \eqref{GREATCONDITIONS} are satisfied if the following conditions hold:
 \begin{equation}\label{GREATCONDITIONS2}
 \left\{\begin{split}
 &\phi k \psi\leq \tau\nu  \\
 &k\mu(\qA_{\supp^c})\leq \frac{1-\nu}{\psi}\\
 &(k-1)\mu(\qA_\supp)\leq 1- \frac{1+\nu}{\psi}\\
 \end{split}\right.
 \end{equation}
 These coherence-based conditions are milder than state-of-the-art, coherence-based CS results, see, e.g., \cite{fuc04,fuc05}, which require $k\leq\frac{1}{2}\left(1+\frac{1}{\mu(A)}\right)$. As a matter of fact, this condition is rather stringent: for example, if $k=2$, one must have $\mu(A)\leq \frac{1}{3}$. In contrast, if we consider $\nu$ very small, we have the condition $k\mu(\qA_{\supp^c})+(k-1)\mu(\qA_\supp)\leq 1$, which allows larger $\mu(\qA)$ provided that $\mu(\qA_\supp)$ is small, which is generally the case when $k\ll n$.
\end{remark}
%
%
\section{Numerical results}\label{sec:nr}
In this section, we present some numerical results. First, we consider a purely static EIV sparse linear regression, with Gaussian $A$. Second, we tackle an EIV ARX system identification problem. In this second experiment, we do not leverage the Toeplitz structure of $A$, defined in \eqref{Toeplitz}. This refinement is left for future work. 

To illustrate the effectiveness of our method, we compare it with the state-of-the-art BPDN$_{\infty}$ \cite{don06infty,fox19acc} and Lasso \cite{tib96,tot11}. In our experiments, we use the alternating direction method of multipliers (ADMM, \cite{boy10}) as convex solver. We choose ADMM because it is easy to implement, even on distributed and parallel architectures, see, e.g., \cite{boy10,mata15,fia18}. However, any convex solver can be used for the purpose.

The performance is evaluated in terms of support recovery, that is, of the identification of the non-zero parameters. This is the substantial problem, as once the support is known, the values of the non-zero parameters ca be assessed, e.g., by  least-squares. Specifically, in our results we show the rates of exact support recovery: we count as  0 a failed support recovery, and as 1 a successful support recovery.
 
For both experiments, we consider $n=100$ and $\Dy=\DA=\Delta$. The proposed results are averaged over 200 random runs.  $\qA$ and $\qy$ are perturbed versions of $A$ and $y$ with maximum error $\Delta$.
\subsection{EIV static sparse linear regression}
In the first experiment, we consider an EIV static sparse linear regression, i.e., the components of $A$ are not correlated. More precisely, we consider   $A\in\R^{m,n}$ with independent components, generated according to a Gaussian  distribution $\mathcal{N}(0,\frac{1}{10})$.   We set $k=10$, and we analyze the performance for different  values of $\Delta$ and $m$. The support of the true $\xtrue$ is generated uniformly at random  while the non-zero entries are uniformly distributed in $[-d,-c]\cup[c,d]$, with $c=\frac{1}{2}$ and $d=1.$
\begin{figure*}
\begin{center}
	\includegraphics[width=0.7\columnwidth]{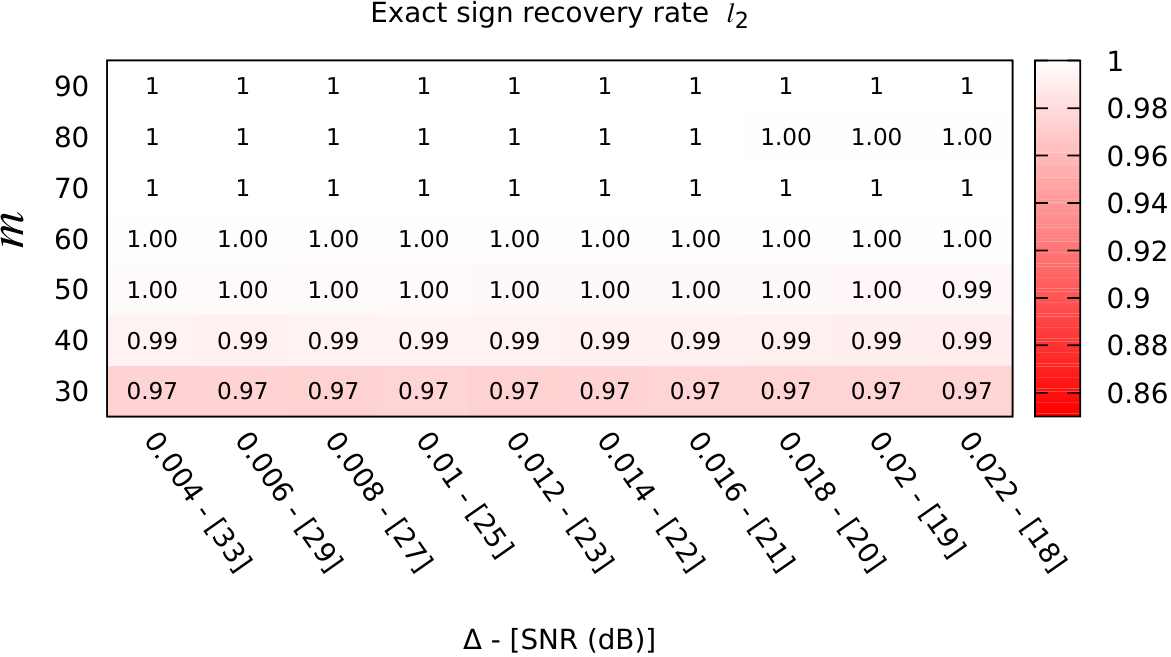}
	\vskip0.4cm
	\includegraphics[width=0.7\columnwidth]{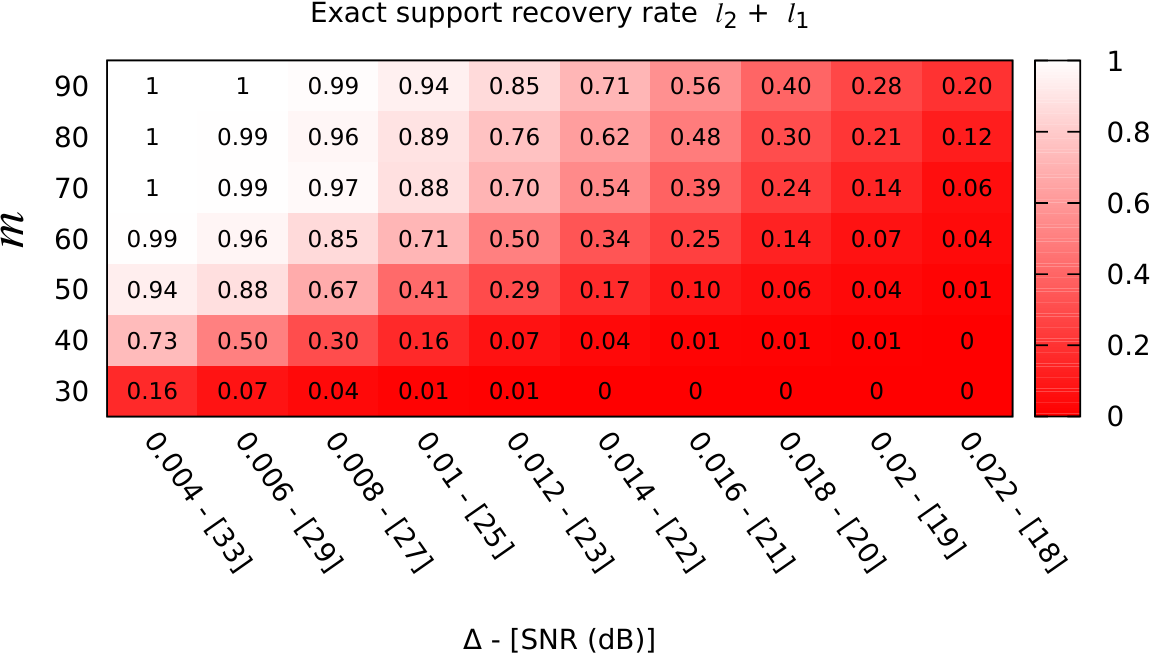}
	\vskip0.4cm
	\includegraphics[width=0.7\columnwidth]{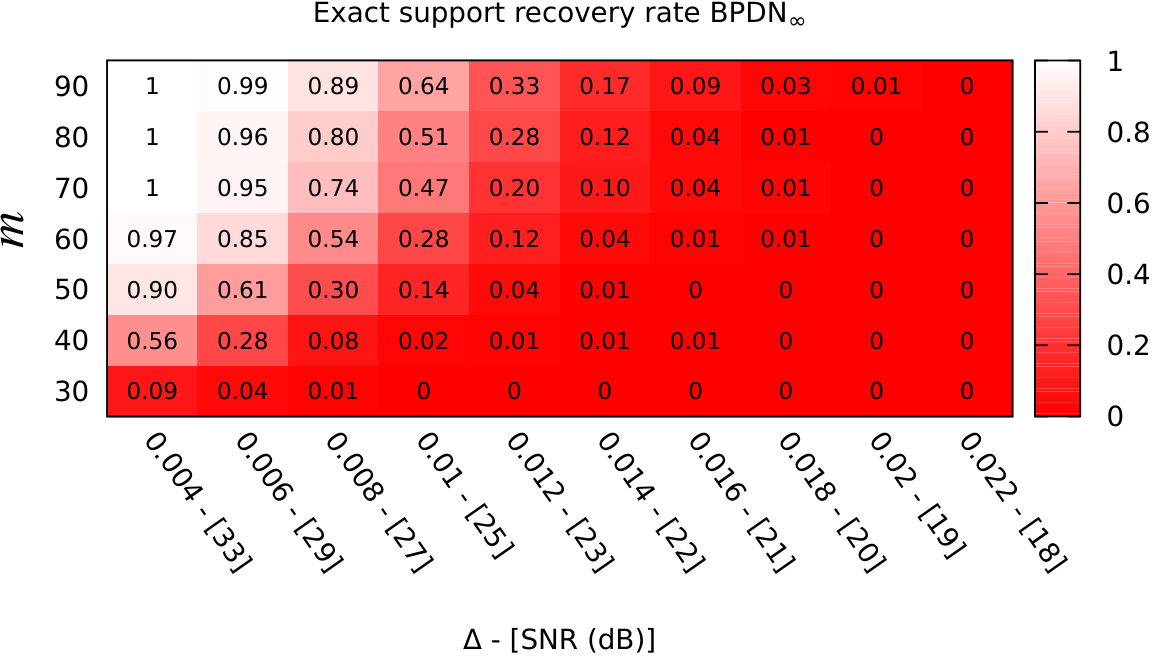}
	\vskip0.4cm
	\includegraphics[width=0.7\columnwidth]{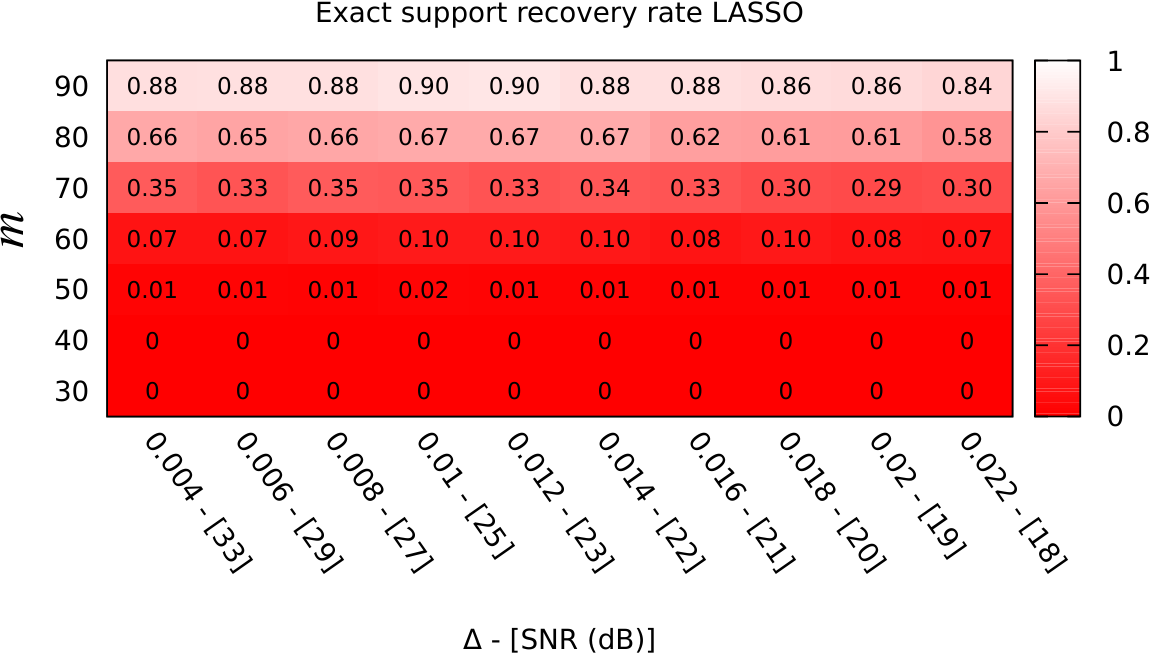}
	\caption{First experiment: EIV linear regression with Gaussian $A$.}
	\label{fig:f}
	\end{center}
\end{figure*}

%
%

In Fig. \ref{fig:f}, we see that the signs $s(\xtrue)$ are recovered with high probability, with no particular sensitivity to noise. Further, in Fig. \ref{fig:f}, we show the exact support recovery rates at different $\Delta=\DA=\Dy$ and $m$. For each $\Delta$, the measured signal-to-noise-ratio  (SNR) is reported as well. The measured SNR (in dB) is defined as $10\log_{10}\frac{\|y\|_2^2+ \|A\|_F^2}{\|\dy\|_2^2+ \|\dA\|_F^2}$, where $\|\cdot\|_F$ is the Frobenius norm. The considered range for $\Delta$ is $[0.004,0.022]$, which corresponds to an SNR varying from 33 to 18 dB.  We can observe that the proposed $\ell_2+\ell_1$ algorithm outperforms both BPDN$_{\infty}$ and Lasso, achieving  more than $90\%$ of success at 25 dB, for a sufficiently large $m$. In particular, Lasso turns out to be definitely ineffective for this fully-perturbed framework. 

\subsection{Identification of ARX system}
In the second experiment, we consider the EIV ARX problem as illustrated in Section \ref{sec:ps}.  Since $A$ has correlated components, the coherence is larger if compared to the first experiments. Therefore, from our theoretical results in Section \ref{sec:ta}, we expect lower performance if compared to Gaussian matrices. This is also observed in Fig. \ref{fig:coherences}. As mentioned above, we do not leverage the Toeplitz structure of $A$. We set $k=10$ and $[c,d]=[0.2, 0.4]$. The values of $c$ and $d$ are smaller than that of previous experiment to guarantee the stability of  ARX system. The input is Gaussian $\mathcal{N}(0,\frac{1}{10})$.
\begin{figure*}
	\centering
	\includegraphics[width=0.7\columnwidth]{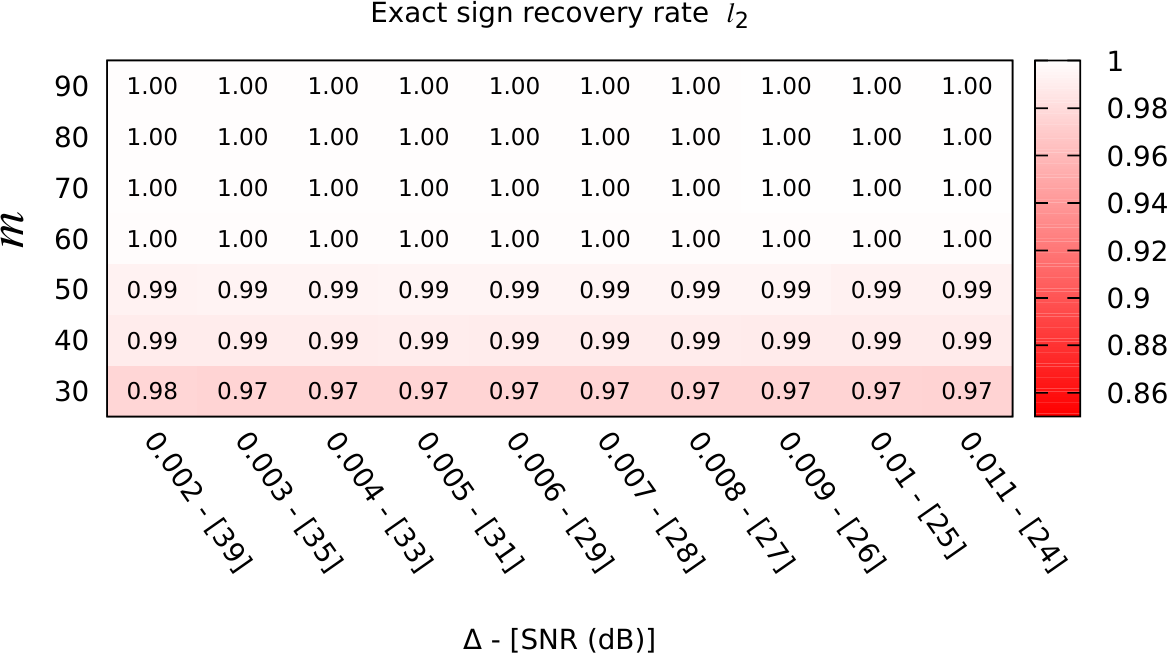}
	\vskip0.4cm
	\includegraphics[width=0.7\columnwidth]{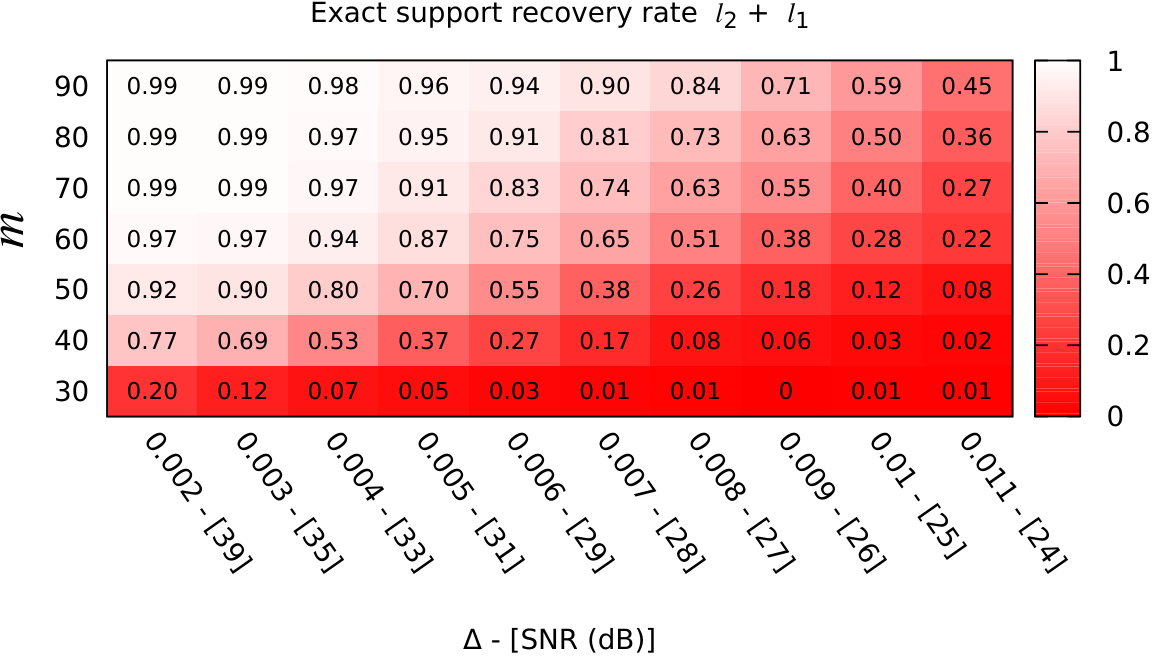}
	\vskip0.4cm
	\includegraphics[width=0.7\columnwidth]{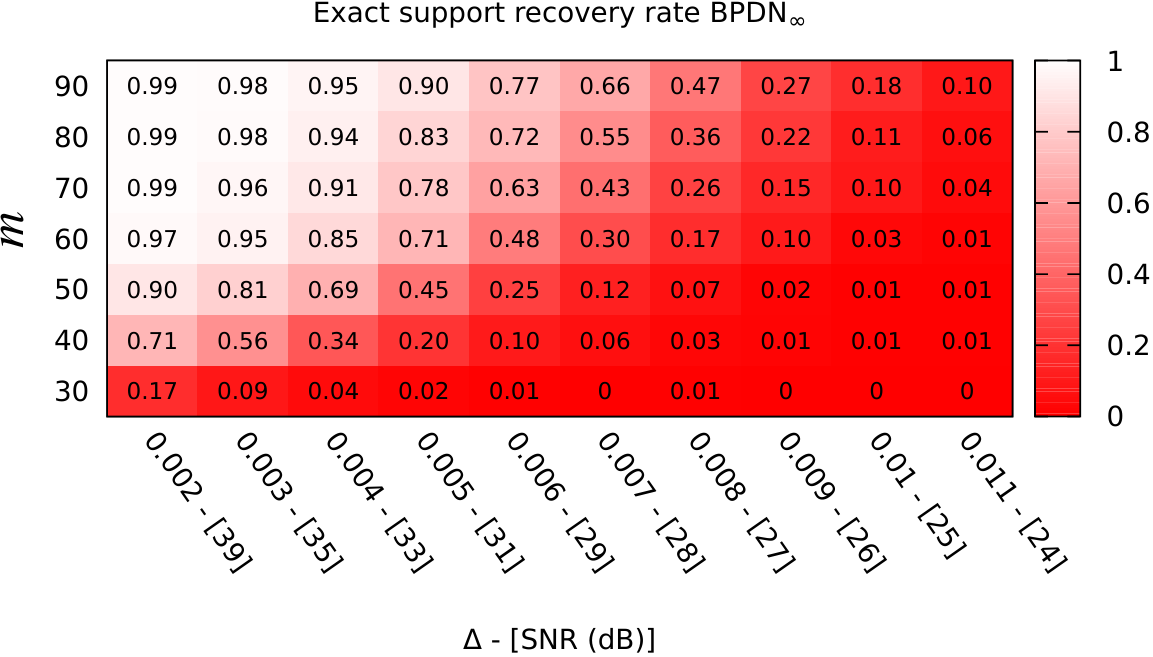}
	\vskip0.4cm
	\includegraphics[width=0.7\columnwidth]{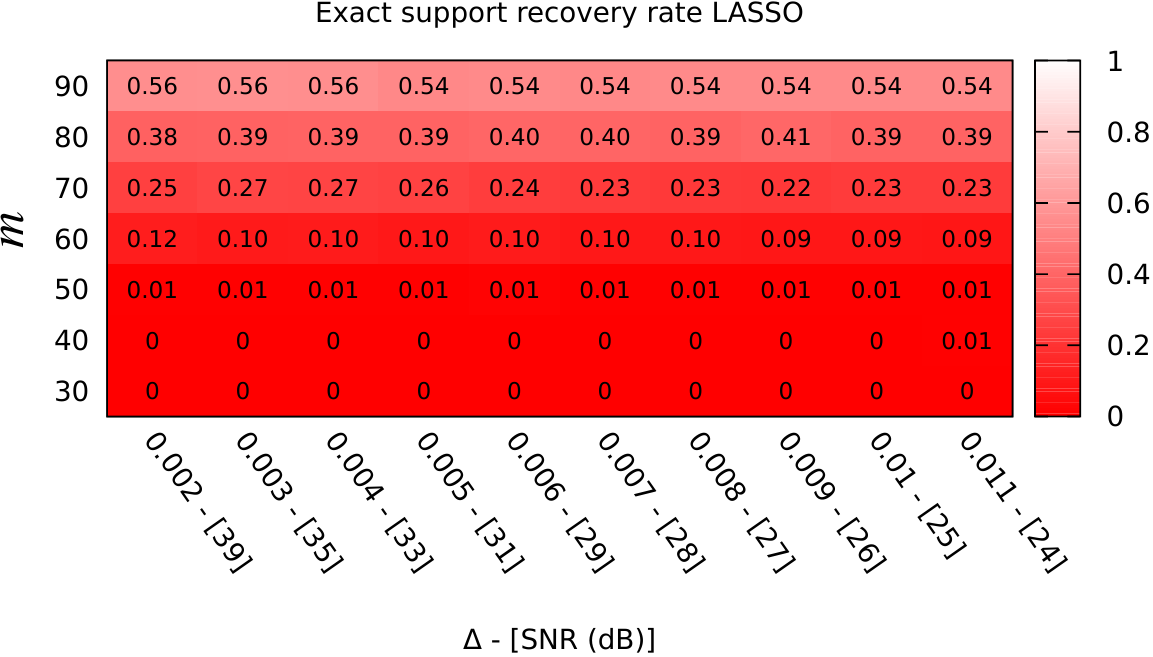}
	\caption{Second experiment: EIV ARX system \eqref{Toeplitz}.}
	\label{fig:g}
\end{figure*}

In Fig. \ref{fig:g}, we can see that the general performance of all the algorithms is slightly degraded if compared to the previous experiment, see Fig. \ref{fig:f}. For this epxeriment, the considered range for $\Delta$ is $[0.002,0.011]$, which corresponds to an SNR varying from 39 to 24 dB. However, the hierarchy is the same: the proposed $\ell_2+\ell_1$ method is always better than BDPN$_{\infty}$ and Lasso. In particular, the $\ell_2+\ell_1$ achieves $90\%$ of success, in terms of support recovery, at 28 dB, for a sufficiently high $m$.
\section{Conclusions}\label{sec:con}
In this paper, we tackle the problem of sparse linear regression from compressed measurements when all the available data are perturbed by noise. The assumption of fully-perturbed data is the most realistic one; nevertheless, the related problem is intrinsically non-convex, thus difficult to solve. In this work, we show that, if perturbations are known to be bounded, an efficient linear programming relaxation is possible. This approach requires to priorly estimate the signs of the solution; for this task, we propose and analyze a Tikhonov approach. The effectiveness of each stage of the proposed strategy is analyzed and sufficient conditions for success are provided. Furthermore, numerical results are proposed to show the performance in practice, on static and dynamic systems. The proposed approach is more effective than the state-of-art methods, in particular for high rate of measurements compression. 
\bibliographystyle{plain}        
\bibliography{refs}    

\begin{thebibliography}{10}

\bibitem{ale08}
A.~Alessandri, M.~Baglietto, and G.~Battistelli.
\newblock Moving-horizon state estimation for nonlinear discrete-time systems:
  New stability results and approximation schemes.
\newblock {\em Automatica}, 44(7):1753--1765, 2008.

\bibitem{bay15}
A.~Bay, D.~Carrera, S.~M. Fosson, P.~Fragneto, M.~Grella, C.~Ravazzi, and
  E.~Magli.
\newblock Block-sparsity-based localization in wireless sensor networks.
\newblock {\em EURASIP J. Wirel. Commun. Netw.}, 2015(182):1--15, 2015.

\bibitem{bem99book}
Al. Bemporad and M.~Morari.
\newblock Robust model predictive control: A survey.
\newblock In A.~Garulli and A.~Tesi, editors, {\em Robustness in identification
  and control}, pages 207--226. Springer London, 1999.

\bibitem{boy10}
S.~Boyd, N.~Parikh, E.~Chu, B.~Peleato, and J.~Eckstein.
\newblock Distributed optimization and statistical learning via the alternating
  direction method of multipliers.
\newblock {\em Found. Trends Mach. Learn.}, 3(1):1 -- 122, 2010.

\bibitem{brubook19}
S.~L. Brunton and J.~Nathan Kutz.
\newblock {\em Data-Driven Science and Engineering: Machine Learning, Dynamical
  Systems, and Control}.
\newblock Cambridge University Press, 2019.

\bibitem{bru16}
S.~L. Brunton, J.~L. Proctor, and J.~N. Kutz.
\newblock Discovering governing equations from data by sparse identification of
  nonlinear dynamical systems.
\newblock {\em Proc. Natl. Acad. Sci.}, 113(15):3932--3937, 2016.

\bibitem{cal16}
M.~Calvo-Fullana, J.~Matamoros, C.~Ant\'{o}n-Haro, and S.~M. Fosson.
\newblock Sparsity-promoting sensor selection with energy harvesting
  constraints.
\newblock In {\em Proc. IEEE Int. Conf. Acoust, Speech Signal Process.
  (ICASSP)}, pages 3766--3770, 2016.

\bibitem{can06}
E.~J. Cand\`es, J.~K. Romberg, and T.~Tao.
\newblock Stable signal recovery from incomplete and inaccurate measurements.
\newblock {\em Commun. Pure Appl. Math.}, 59(8):1207--1223, 2006.

\bibitem{car13}
A.~Y. Carmi.
\newblock Compressive system identification: Sequential methods and entropy
  bounds.
\newblock {\em Digital Signal Process.}, 23(3):751--770, 2013.

\bibitem{cer93}
V.~Cerone.
\newblock Feasible parameter set for linear models with bounded errors in all
  variables.
\newblock {\em Automatica}, 29(6):1551 -- 1555, 1993.

\bibitem{fox19acc}
V.~Cerone, S.~M. Fosson, and D.~Regruto.
\newblock A linear programming approach to sparse linear regression with
  quantized data.
\newblock In {\em Proc. American Control Conf. (ACC)}, pages 2990--2995, 2019.

\bibitem{pig12}
V.~{Cerone}, D.~{Piga}, and D.~{Regruto}.
\newblock Set-membership error-in-variables identification through convex
  relaxation techniques.
\newblock {\em IEEE Trans. Autom. Control}, 57(2):517--522, 2012.

\bibitem{cer17}
V.~{Cerone}, D.~{Regruto}, and M.~{Abuabiah}.
\newblock Direct data-driven control design through set-membership
  errors-in-variables identification techniques.
\newblock In {\em Proc. Amer. Control Conf. (ACC)}, pages 388--393, 2017.

\bibitem{cic19}
V.~{Ciccone}, A.~{Ferrante}, and M.~{Zorzi}.
\newblock Factor models with real data: A robust estimation of the number of
  factors.
\newblock {\em IEEE Trans. Autom. Control}, 64(6):2412--2425, 2019.

\bibitem{dat19}
S.~Datta and J.~Oldroyd.
\newblock Low coherence unit norm tight frames.
\newblock {\em Linear and Multilinear Algebra}, 67(6):1174--1189, 2019.

\bibitem{don06infty}
D.~L. Donoho and M.~Elad.
\newblock On the stability of the basis pursuit in the presence of noise.
\newblock {\em Signal Processing}, 86(3):511 -- 532, 2006.

\bibitem{don01}
D.~L. {Donoho} and X.~{Huo}.
\newblock Uncertainty principles and ideal atomic decomposition.
\newblock {\em IEEE Trans. Inf. Theory}, 47(7):2845--2862, 2001.

\bibitem{ela07}
M.~{Elad}.
\newblock Optimized projections for compressed sensing.
\newblock {\em IEEE Trans. Signal Process.}, 55(12):5695--5702, 2007.

\bibitem{fat18}
S.~{Fattahi} and S.~{Sojoudi}.
\newblock Data-driven sparse system identification.
\newblock In {\em Proc. Allerton Conf. Commun. Control Comput.}, pages
  462--469, 2018.

\bibitem{fen09}
C.~Feng, S.~Valaee, and Z.~Tan.
\newblock Multiple target localization using compressive sensing.
\newblock In {\em Proc. IEEE Global Telecommun. Conf. (GLOBECOM)}, pages 1--6,
  2009.

\bibitem{fia18}
A.~Fiandrotti, S.~M. Fosson, C.~Ravazzi, and E.~Magli.
\newblock {GPU}-accelerated algorithms for compressed signals recovery with
  application to astronomical imagery deblurring.
\newblock {\em Int. J. Remote Sens.}, 39(7):2043--2065, 2018.

\bibitem{fou13}
S.~Foucart and H.~Rauhut.
\newblock {\em A Mathematical Introduction to Compressive Sensing}.
\newblock Springer, New York, 2013.

\bibitem{fra17}
G.~Fracastoro, S.~M. Fosson, and E.~Magli.
\newblock Steerable discrete cosine transform.
\newblock {\em IEEE Trans. Image Process.}, 26(1):303--314, 2017.

\bibitem{fuc04}
J.~J. Fuchs.
\newblock On sparse representations in arbitrary redundant bases.
\newblock {\em IEEE Trans. Inf. Theory}, 50(6):1341--1344, 2004.

\bibitem{fuc05}
J.~J. Fuchs.
\newblock Recovery of exact sparse representations in the presence of bounded
  noise.
\newblock {\em IEEE Trans. Inf. Theory}, 51(10):3601--3608, 2005.

\bibitem{voy16}
L.~Gallana, F.~Fraternale, M.~Iovieno, S:~M. Fosson, E.~Magli, M.~Opher, J.~D.
  Richardson, and D.~Tordella.
\newblock Voyager 2 solar plasma and magnetic field spectral analysis for
  intermediate data sparsity.
\newblock {\em J. Geophys. Res. Space Phys.}, 121(5):3905--3919, 2016.

\bibitem{gol73}
G.~H. Golub.
\newblock Some modified matrix eigenvalue problems.
\newblock {\em SIAM Review}, 15(2), 1973.

\bibitem{gol80}
G.~H. Golub and C.~F. van Loan.
\newblock An analysis of the total least squares problem.
\newblock {\em SIAM J. Numer. Anal.}, 17(6):883--893, 1980.

\bibitem{gu09}
Y.~Gu, J.~Jin, and S.~Mei.
\newblock $\ell_{0}$ norm constraint {LMS} algorithm for sparse system
  identification.
\newblock {\em IEEE Signal Process. Lett.}, 16(9):774--777, 2009.

\bibitem{her10}
M.~A. Herman and T.~Strohmer.
\newblock General deviants: An analysis of perturbations in compressed sensing.
\newblock {\em IEEE J. Sel. Top. Sign. Proces.}, 4(2):342--349, 2010.

\bibitem{las01}
J.~Lasserre.
\newblock Global optimization with polynomials and the problem of moments.
\newblock {\em SIAM J. Optim.}, 11(3):796--817, 2001.

\bibitem{lasbook}
J.-B. Lasserre.
\newblock {\em An introduction to polynomial and semi-algebraic optimization}.
\newblock Cambridge University Press, 2015.

\bibitem{lue16}
D.~G. Luenberger and Y.~Ye.
\newblock {\em Linear and Nonlinear Programming}.
\newblock Springer International Publishing Switzerland, 4th edition, 2016.

\bibitem{mar07}
I.~Markovsky and S.~Van Huffel.
\newblock Overview of total least-squares methods.
\newblock {\em Signal Process.}, 87(10):2283--2302, 2007.

\bibitem{mata15}
J.~Matamoros, S.~M. Fosson, E.~Magli, and C.~Ant\'{o}n-Haro.
\newblock Distributed {ADMM} for in-network reconstruction of sparse signals
  with innovations.
\newblock {\em IEEE Trans. Signal Inf. Process. Netw.}, 1(4):225--234, 2015.

\bibitem{ric05}
A.~{Richards} and J.~{How}.
\newblock Robust model predictive control with imperfect information.
\newblock In {\em Proc. Amer. Control Conf.}, pages 268--273, 2005.

\bibitem{roj14}
C.~R. {Rojas}, R.~{T\'{o}th}, and H.~{Hjalmarsson}.
\newblock Sparse estimation of polynomial and rational dynamical models.
\newblock {\em IEEE Trans. Autom. Control}, 59(11):2962--2977, 2014.

\bibitem{san11}
B.~M. Sanandaji, T.~L. Vincent, M.~B. Wakin, and R.~T\'{o}th.
\newblock Compressive system identification of lti and ltv arx models.
\newblock In {\em Proc. IEEE Conf. Decision Control (CDC)}, pages 783--790,
  2011.

\bibitem{sod07}
T.~S\"oderstr\"om.
\newblock Errors-in-variables methods in system identification.
\newblock {\em Automatica}, 43(6):939 -- 958, 2007.

\bibitem{sod18}
T.~S\"oderstr\"om.
\newblock {\em Errors-in-Variables Methods in System Identification}.
\newblock Springer, 2018.

\bibitem{tib96}
R.~Tibshirani.
\newblock Regression shrinkage and selection via the {L}asso.
\newblock {\em J. Royal. Statist. Soc. B}, 58:267--288, 1996.

\bibitem{tot11}
R.~T\'{o}th, B.~M. Sanandaji, K.~Poolla, and T.~L. Vincent.
\newblock Compressive system identification in the linear time-invariant
  framework.
\newblock In {\em Proc. IEEE Conf. Decision Control (CDC)}, pages 783--790,
  2011.

\bibitem{tsi14}
E.~V. {Tsiligianni}, L.~P. {Kondi}, and A.~K. {Katsaggelos}.
\newblock Construction of incoherent unit norm tight frames with application to
  compressed sensing.
\newblock {\em IEEE Trans. Inf. Theory}, 60(4):2319--2330, 2014.

\bibitem{yan12}
Z.~Yang, C.~Zhang, and L.~Xie.
\newblock Robustly stable signal recovery in compressed sensing with structured
  matrix perturbation.
\newblock {\em IEEE Trans. Signal Process.}, 60(9):4658--4671, 2012.

\bibitem{zhu11}
H.~Zhu, G.~Leus, and G.~B. Giannakis.
\newblock Sparsity-cognizant total least-squares for perturbed compressive
  sampling.
\newblock {\em IEEE Trans. Signal Process.}, 59(5):2002--2016, 2011.

\end{thebibliography}
\end{document}